\documentclass[12pt]{article}
\usepackage{latexsym,amsmath,amssymb,amsfonts}
\usepackage{graphicx}

\newcommand{\vep}{\mbox{$\varepsilon$}}
\newcommand{\pa}{\mbox{$\partial$}}

\newcommand{\ds}{\displaystyle}

\newcommand{\vslash}{\mbox{\,\rule[-0.37cm]{0.020cm}{1.0cm}\,}}

\newcommand{\YY}{\mbox{$\widetilde{Y}$}}

\newtheorem{eg}{Example}[section]
\newtheorem{thm}{Theorem}[section]
\numberwithin{equation}{section}
\begin{document}
	
	\large
	
	\begin{flushleft}
		{\large \bf Some remarks on the solution of linearisable second-order ordinary differential equations via point transformations}\\[3 mm]		
		
		{ Winter Sinkala\\[3 mm]}
		Department of Mathematical Sciences and Computing\\ Faculty of Natural Sciences, 
		Walter Sisulu University\\ Private Bag X1, Mthatha 5117, Republic of
		South Africa\\ (E-mail address: wsinkala@wsu.ac.za).
	\end{flushleft}

	\vspace{-5mm}

	\begin{abstract}
		\noindent Transformations of differential equations to other equivalent equations play a central role in many routines for solving intricate equations. A class of differential equations that are particularly amenable to solution techniques based on such transformations is the class of linearisable second-order ordinary differential equations (ODEs). There are various characterisations of such ODEs. We exploit a particular characterisation and the expanded Lie group method to construct a generic solution for all linearisable second-order ODEs. The general solution of any given equation from this class is then easily obtainable from the generic solution through a point transformation constructed using only two suitably chosen symmetries of the equation. We illustrate the approach with three examples.
	\end{abstract}
	
	\smallskip
	
	{\small KEY WORDS: Lie symmetry, Nonlinear second-order ODE, Expanded Lie group, Linearising transformation.}
	
	
	\section{Introduction}
	
	Nonlinear ODEs arise in many different contexts, including as mathematical models of real-world phenomena. Analytical solutions of many such equations are often hard to find, which is why a whole range of methods have been proposed for investigating different types of nonlinear ODEs. These methods include Painlev\'{e} singularity analysis, Lie symmetry analysis, Darboux method, and the Jacobi last multiplier method (see \cite{Senthilvelan} and the references therein). The Lie symmetry method, which is based on the invariance of a differential equation under a continuous group of point transformations, is widely used. Given a differential equation, Lie point symmetries of the equation can be used to perform many things on the equation including constructing explicit transformations that reduce the equation to a simpler one, when this is possible. When the simpler ``target" equation is a linear equation the problem is called the \emph{linearisation problem}. The pioneering work on this, with regard to second-order ODEs, is attributed to Sophus Lie (see \cite{SAFDAR11} and the references there in). Lie proved that to be linearisable, a second-order ODE must be at most cubically semi-linear and the coefficients in it must satisfy an overdetermined system of conditions \cite{Mahomed07b, Schwarz}. A considerable amount of research has since been conducted on the linearisation problem
	\cite{Mahomed07b,MUSTAFA13,Safdar,SAFDAR11,Nakpim,Safdar,Nakpim16} leading, in particular, to a variety of ways of characterising linearisable second-order ODEs (see, e.g. \cite{Mahomed07b} [Theorem 8]). 
	
	We remark here that the study of differential equations in general is a diverse field and has attracted intensive research over the years. Important recent advances in the area has included development of new solution methods to solve various classes of nonlinear partial differential equations (see \cite{Mahmoud18, Mahmoud18b, Mahmoud20, Mahmoud20b, Mahmoud20c} and the references therein). 
	
In respect of characterising linearisable second-order ODEs, a scaled down version of Theorem~8 of \cite{Mahomed07b}  is presented here.
	\begin{thm}\label{mhthm}
		Let us consider a scalar second-order ODE in the form
		\begin{equation}\label{fm17}
		y^{\prime\prime} = f(x,y,y^\prime).
		\end{equation}
		The following statements are equivalent \cite{Mahomed07b}:
		\begin{enumerate}
			\item A scalar second-order Eq. (\ref{fm17}) is linearizable via a point transformation.
			\item Eq. (\ref{fm17}) has the maximum eight-dimensional symmetry Lie algebra.
			\item\label{bullet3} Eq. (\ref{fm17}) has the cubic in derivative form
			\begin{equation}\label{fm19}
			y^{\prime\prime} = A(x, y) {y^\prime}^3+B(x, y) {y^\prime}^2+C(x, y) {y^\prime}+D(x, y),
			\end{equation}	
			with the coefficients $A$ to $D$ satisfying
			the two invariant conditions
			\begin{equation}
			\begin{array}{rr}
			3A_{xx} +3A_xC-3A_yD+3AC_x +C_{yy}-6AD_y+BC_y-2BB_x - 2B_{xy} = 0 &  \\
			6A_x D-3B_yD+3AD_x +B_{xx} - 2C_{xy}-3BD_y+3D_{yy}+2C C_y-CB_x = 0. &
			\end{array}
			\end{equation}
			\item\label{bullet4} Eq. (\ref{fm17}) has two non-commuting symmetries $X_1$, $X_2$, in a suitable basis with
			\begin{equation}\label{e4w}
			[X_1, X_2]=X_1, \quad X_1\ne\rho(x,y)X_2
			\end{equation}
			for any non-constant function $\rho$ such that a point change of variables $X =X(x, y), Y =Y (x, y)$ which brings $X_1$ and $X_2$ to their canonical form
			\begin{equation}\label{cform}
			X_1 = \frac{\partial}{\partial Y}, \quad X_2 = X \frac{\partial}{\partial X} +Y\frac{\partial}{\partial Y}
			\end{equation}
			reduces Eq. (\ref{fm17}) to
			\begin{equation}\label{e31}
			X Y^{\prime\prime} = \frac{b}{3 a} + \frac{b^3}{27 a^2} + \left(1 + \frac{b^2}{3 a}\right) Y^\prime + b {Y^\prime}^2 + a {Y^\prime}^3
			\end{equation}
			where $a(\ne 0)$ and $b$ are constants.
		\end{enumerate}
	\end{thm}
	In this paper, we take advantage of the representation (\ref{e31}) to construct a proxy solution for all linearisable second-order ODEs. We start by using the expanded Lie group method to simplify Eq. (\ref{e31}) significantly. We then construct two point transformations and use them to recover the general solution of Eq. (\ref{e31}) from the solution of the simplest second-order ODE, the free particle equation $Y^{\prime\prime}=0$.
	The solution of any given linearisable second-order ODE can now be recovered from the solution of Eq. (\ref{e31}) in a routine fashion via a point transformation constructed from only two symmetries of the equation. We illustrate the solution routine through three examples.
	
	The rest of the paper is organised as follows. In Section \ref{sec2} we use the expanded Lie group method to construct a point transformation that reduces Eq. (\ref{e31}) to a simpler equation. In Section \ref{sec3} we perform further reduction and map of Eq. (\ref{e31}) to the free particle equation. We subsequently deduce a generic solution of all linearisable second-order ODEs.  In Section \ref{sec4} we illustrate through three examples how the general solution of any linearisable second-order ODE may be deduced  from the constructed generic solution. We give concluding remarks in Section \ref{sec5}.

	\section{Reduction of Eq. (\ref{e31}) via the expanded Lie group method}\label{sec2}
	The Lie symmetry method for studying differential equations, initiated by Sophus Lie in the latter part of the nineteenth century, is based on continuous groups of transformations that map solutions of a given differential equation into other solutions of the same equation. The method extends and harmonises various specialised methods for solving ODEs. There is extensive literature on the Lie symmetry method, to which we refer the interested reader (see, for example \cite{Schwarz,Baumann2000,BlumanKumei,Cantwell,Hydon2000,OvsiannikovLV82,olver}).
	
	When we consider a continuous group of transformations acting on the expanded space of variables, which includes
	the equation parameters in addition to independent and dependent variables, we obtain an expanded Lie group transformation of the equation \cite{burde}.
	Such a group of transformations represents a particular case of the equivalence group that preserves the class of equations under study.

	Let us take the parameter $b$ in Eq. (\ref{e31}) as a second independent variable. Now consider a one-parameter ($\vep$) Lie group of point transformations in $(X, Y, b)$:
	\begin{equation}\label{InvPaperLgp}
	\widetilde{X} = f(X, Y, b,\vep), \quad
	\YY = g(X, Y, b,\vep), \quad
	\widetilde{b} = h(b,\vep),
	\end{equation}
	for some functions $f$, $g$ and $h$, with an infinitesimal generator of the form
	\begin{equation}\label{DeltaOperator}
	X = \xi(X, Y, b)\,\pa_X + \tau(X, Y, b)\,\pa_Y + \beta(b)\,\pa_b,
	\end{equation}
	which leaves Eq. (\ref{e31}) invariant. Eq. (\ref{e31}) admits (\ref{InvPaperLgp}) if  $X$ satisfies the invariance requirement
	\begin{equation}\label{InvCondBond}
	X^{[2]}\left\{\frac{b}{3 a} + \frac{b^3}{27 a^2} + \left(1 + \frac{b^2}{3 a}\right) Y^\prime + b {Y^\prime}^2 + a {Y^\prime}^3-X Y^{\prime\prime}\right\}\vslash_{(\ref{e31})} = 0,
	\end{equation}
	where $X^{[2]}$ is the second extension of $X$.
	We obtain, as a particular solution of (\ref{InvCondBond}), that
	\begin{eqnarray}\label{}
	\xi &=& \frac{1 + X^2}{2\,X},\\
	\tau &=& Y\bigg[\frac{1 + X^2}{4\,X} + \frac{1}{2}\left(1-\frac{1 + X^2}{2\,X^2}\right)\bigg]-\frac{b}{6\,a\,X} - \frac{X}{3\,a}\\
   \beta &=& 1.
	\end{eqnarray} The corresponding one-parameter ($\vep$) Lie group of transformations is
	\begin{eqnarray}
	\widetilde{X}  &=& \sqrt{e^{\vep }\,\left( 1 + X^2 \right) - 1}\label{e5}\\
  	\YY  &=& \frac{e^{\vep/2}(3\,a\,Y + b\,{X}) -
  	\left( b + \vep  \right) {\sqrt{e^{\vep }\left( 1 + X^2 \right) - 1}} }{3\,a}\label{e6}\\
	\widetilde{b} &=& b + \vep,\label{e7}
	\end{eqnarray}
	so that under this transformation Eq. (\ref{e31}) necessarily becomes
	\begin{equation}\label{e31b}
	\widetilde{X} \YY^{\prime\prime} = \frac{\widetilde{b}}{3 a} + \frac{\widetilde{b}^3}{27 a^2} + \left(1 + \frac{\widetilde{b}^2}{3 a}\right) \YY^\prime + \widetilde{b} {\YY^\prime}^2 + a {\YY^\prime}^3.
	\end{equation}
	If we now set $\vep = -b$ in the transformation (\ref{e5})--(\ref{e7}), $\widetilde{b}$ equals zero and Eq. (\ref{e31b}) reduces to
	\begin{equation}\label{e31c}
	\widetilde{X} \YY^{\prime\prime} = \YY^\prime +  a {\YY^\prime}^3.
	\end{equation}
	Clearly Eq. (\ref{e31c}) is mapped back to Eq. (\ref{e31}) via the point transformation (\ref{e5})--(\ref{e7}) with $\vep = -b$, i.e.
	\begin{equation}\label{e31d}
	\widetilde{X}  = \sqrt{e^{-b}\,\left( 1 + X^2 \right) - 1}, \quad
	\YY =  e^{-b/2}\left(Y + \frac{b}{3\,a}\,{X}\right).
	\end{equation}
	This means that a solution of Eq. (\ref{e31c}) is mapped to a solution of Eq. (\ref{e31}) via the transformation (\ref{e31d}).

	\section{Reduction of  (\ref{e31c}) to the free particle equation}\label{sec3}
	We seek an invertible point transformation between Eq. (\ref{e31c}), which we restate here in the original variables $X$ and $Y$,
	\begin{equation}\label{e31e}
	X Y^{\prime\prime} = Y^\prime +  a {Y^\prime}^3,
	\end{equation}
	and the free particle equation
	\begin{equation}\label{fpe}
	W^{\prime\prime} = 0, \quad W = W(Z).
	\end{equation}
	We exploit the equivalence of the symmetry Lie algebras of the two equations to construct the point transformation.
	
	Lie point symmetries admitted by equations (\ref{e31e})  and (\ref{fpe}) are
	\begin{equation}\label{21ww}
	 \begin{array}{l}
	\Phi_1 = \frac{1}{X}\,\pa_X, \quad \Phi_2 = (-2\,a^2\,Y^3/X)\pa_X + (3\,a\,Y^2 + X^2)\pa_Y\\[1ex]
	\Phi_3 = \left(a\,Y^2/X+ X\right)\pa_X, \quad \Phi_4 = (a\,Y^2/X)\pa_X - Y \pa_Y\\[1ex]
	\Phi_5 = (X^4-a^2\,Y^4/2X)\pa_X + Y\,\left( a\,Y^2 + X^2 \right)\pa_Y, \quad \Phi_6 = (Y/X)\pa_X\\[1ex]
	\Phi_7 = (2\,a\,Y/X)\pa_X - \pa_Y, \quad \Phi_8 = (Y\,X - a\,Y^3/X)\pa_X + 2\,Y^2\pa_Y\\
	\end{array}
	\end{equation}
	and
	\begin{equation}\label{22ww}
	\begin{array}{l}
	\Omega_1 = \partial_Z, \quad \Omega_2 = \partial_W, \quad \Omega_3 = Z\partial_Z, \quad \Omega_4 = Z\partial_W, \quad \Omega_5 = W\partial_W  \\[1ex]
	\Omega_6 = W\partial_Z, \quad \Omega_7 = Z W\partial_Z + W^2\partial_W, \quad \Omega_8 = Z^2\partial_Z-Z W\partial_W,
	\end{array}
	\end{equation}
	respectively. The Lie algebras arising from (\ref{21ww}) and (\ref{22ww}) are equivalent, which means an isomorphism leading to the
	same commutator table for the two Lie algebras can be found. We may rearrange the operators (\ref{21ww}) to form a new basis $\{\Gamma_1, \ldots, \Gamma_8\}$ such that
	\[[\Gamma_i ,\Gamma_j]=C_{i j}^k \Gamma_k\quad  \mbox{and}\quad [\Omega_i ,\Omega_j]=C_{i j}^k \Omega_k,\]
	with the same structure constants $\{C_{i j}^k\}$. The following rearrangements provide the desired new basis for the symmetry Lie algebra arising from (\ref{21ww}):
\begin{equation}\label{}
\begin{array}{l}
\ds\Gamma_1 = {\gamma} \Phi_6, \quad
\Gamma_2 = {\gamma}\,{\lambda} \Phi_8, \quad
\Gamma_3 = \frac{1}{2} \Phi_3, \quad
\Gamma_4 = {\lambda} \Phi_4,\\[1ex]
\ds\Gamma_5 = -\frac{{\Phi_3}}{2} + {\Phi_4} +
{\gamma}\,{\delta}\,{\Phi_8}, \quad
\Gamma_6 = -\frac{{\Phi_1}}{2\,{\lambda}} +
\frac{{\gamma}\,{\delta}\,
	{\Phi_6}}{{\lambda}},\\[1ex]
\ds\Gamma_7 = -\frac{{\delta}\,{\Phi_3}}
{2\,{\lambda}} +
\frac{2\,{\delta}\,{\Phi_4}}
{{\lambda}} +
\frac{a\,{\Phi_6}}
{2\,{\gamma}\,{\lambda}} -
\frac{{\Phi_7}}
{2\,{\gamma}\,{\lambda}} +
\frac{{\gamma}\,{{\delta}}^2\,
	{\Phi_8}}{{\lambda}},\\[1ex]
\ds\Gamma_8 = -\frac{{\Phi_2}}{2\,{\gamma}} +
{\delta}\,{\Phi_5} +
\frac{a\,{\Phi_8}}{2\,{\gamma}},
\end{array}
\end{equation}

where $\gamma$, $\lambda$ and $\delta$ are arbtrary constants with $\lambda \gamma\ne 0$.
	We now seek a point transformation that maps Eq. (\ref{e31e}) to Eq. (\ref{fpe}) of the form
	\begin{equation}\label{3.5}
	Z = \alpha (X,Y), \quad W = \beta(X,Y)
	\end{equation}
	for functions $\alpha$ and $\beta$. According to \cite[Chapter 6]{BlumanKumei} the functions $\alpha$ and $\beta$ must be such that the conditions
	\begin{equation}\label{24ww}
	\Omega_i Z = \Gamma_i\, \alpha (X,Y), \quad \Omega_i W = \Gamma_i\, \beta(X,Y), \quad i=1,\ldots,8
	\end{equation}
	are satisfied. The equations in (\ref{24ww}) translate into an overdetermined system of sixteen elementary linear PDEs that define the functions $\alpha$ and $\beta$. The equations are easily solved and we obtain
	\begin{equation}\label{24w}
	Z = \alpha (X,Y) = \frac{a\,Y^2 + X^2}{2\,{\gamma}\,Y}, \quad W = \beta (X,Y) = \frac{{\delta}}{{\lambda}} -
	\frac{1}{2\,{\gamma}\,{\lambda}\,Y}.
	\end{equation}
	The solution of (\ref{e31e}) is now recovered from the solution of (\ref{fpe}),
	\begin{equation}\label{25w}
	W = m Z + c,
	\end{equation}
	through the point transformation (\ref{24w}). We obtain
	\begin{equation}\label{26w}
	 \frac{{\delta}}{{\lambda}} -
	\frac{1}{2\,{\gamma}\,{\lambda}\,Y}   =  m \left(\frac{a\,Y^2 + X^2}{2\,{\gamma}\,Y} \right)  + c,
	\end{equation}
	or, equivalently,
	\begin{equation}\label{27w}
	a\,Y^2 + A\,Y + X^2 + B  = 0,
	\end{equation}
	where \[A = \frac{2\,\gamma}{m}\left(c - \frac{\delta}{\lambda}\right)  \quad \mbox{and} \quad B = \frac{1}{{\lambda}\,m}.\]
	The solution of the generic Eq. (\ref{e31}) now follows from the point transformation (\ref{e31d}) and Eq. (\ref{27w}), stated in terms of the variables $\widetilde{X}$ and $\YY$ of Eq. (\ref{e31c}). It is
	\begin{equation}\label{28w}
	a\, {Y}^2 +  \frac{2\, b} {3}\, X \,
	Y + \left (1 + \frac{b^2} {9\, a} \right)\,
	X^2 + J_2\left( Y + \frac{b} {3\, a}\, X \right) + J_1 = 0,
	\end {equation}
	where \[J_1 = e^b (B-1)+1, \quad J_2 = A e^{b/2}.\]
	Eq. (\ref{28w}) represents a convenient solution of every linearisable second-order ordinary differential equation in that the general solution of every such equation is realisable from Eq. (\ref{28w}) via a point transformation constructed using only two suitably chosen symmetries of the equation according to Theorem~\ref{mhthm}.
	\section{Illustrative examples}\label{sec4}
	\begin{eg} 
	{\rm
		The ODE
		\begin{equation}\label{29ww}
		y^{\prime\prime} + 3 y y^\prime + y^3 = 0
		\end{equation}
		called the modified Emdem equation, which arises in a vriety of contexts \cite{Senthilvelan}, admits the maximal 8 symmetries. Among the admitted symmetries are two noncomuting symmetry generators
		\begin{equation}\label{31ww}
		X_1 = \pa_x \quad\mbox{and}\quad  X_2 = x \pa_x - y \pa_y
		\end{equation}
		that satisfy condition (\ref{e4w}) of Theorem~\ref{mhthm}. The symmetries (\ref{31ww}) are reduced to their canonical form (\ref{cform}) via the point transformation
		\begin{equation}\label{32ww}
		X= k_1/y, \quad Y = x + k_2/y,
		\end{equation}
		where $k_1(\ne 0)$ and $k_2$ are arbitrary constants (see also \cite{Mahomed07b}).
		Under this point transformation Eq. (\ref{29ww}) is reduced to Eq. (\ref{e31}) with $a=-k_1^2$ and $b=3 k_1 (k_2+1)$, and the solution in (\ref{28w}) is transformed into
\begin{equation}\label{}
    y= \frac{2 k_1^2 x - J_2}{k_1^2 x^2 - J_2 x - J_1},
\end{equation}
		which is the desired general solution of (\ref{29ww}).
	}
\end{eg}
\begin{eg} 
	{\rm
		Consider ODE No. 6.180 of Kamke \cite{Kamke}:
		\begin{equation}\label{34ww}
		y^{\prime\prime}(y-1)^2 x^2 -  2 {y^\prime}^2x^2 - y^\prime(y-1) x - 2y(y-1)^2 = 0.
		\end{equation}
		This equation satisfies the conditions in Item~\ref {bullet3} of Theorem~\ref{mhthm}
		and therefore admits the maximal $8$ Lie point symmetries, two of which are:
		\begin{equation}\label{35ww}
		X_1 = \left(\frac{y-2}{y-1}\right) \pa_x + \frac{2 y}{x}\,\pa_y \quad\mbox{and}\quad X_2 = x \pa_x.
		\end{equation}
		These infinitesimal symmetries satisfy condition (\ref{e4w}) from Theorem~\ref{mhthm} and are reduceable to their canonical forms (\ref{cform}) via the point transformation
		\begin{equation}\label{36ww}
		X= \frac{k_2 x \sqrt{1-y}}{y}, \quad Y = \left(\frac{k_1 \sqrt{1-y}}{y} + 1 -\frac{1}{y}\right)x,
		\end{equation}
		where $k_1$ and $k_2(\ne 0)$ are arbitrary constants.
		Under this transformation Eq. (\ref{34ww}) is reduced  to  Eq. (\ref{e31}), with $a=-k_2^2$ and $b=3 k_1 k_2$. The general solution of (\ref{34ww}) now follows easily from Eq. (\ref{28w}) written in terms of $x$ and $y$ via the point transformation (\ref{36ww}), with $a$ and $b$ set to $-k_2^2$ and $3 k_1 k_2$, respectively. We obtain
\begin{equation}\label{}
	y = \frac{x \left(k_2^2 x - J_2\right)}{k_2^2 x^2 - J_2 x - J_1}.
\end{equation}	
	
	}
\end{eg}
\begin{eg} 
	{\rm
		Consider the linearisable equation \cite[p. 100]{Hydon2000}
		\begin{equation}\label{38ww}
		x y^{\prime\prime} = 2\left(y^{\prime } \right)^2 + y^{\prime } + y^2.
		\end{equation}
		The symmetries
		\begin{equation}\label{39ww}
		X_1 = \left(1-\frac{1}{x y} \right)\pa_x + \frac{y}{x}\,\pa_y  \quad\mbox{and}\quad X_2 = \left(x^2-\frac{x}{y} \right)\pa_x + \left(2-3 x y \right)\pa_y
		\end{equation}
		are admitted by (\ref{38ww}) and satisfy the condition (\ref{e4w}). The point transformation that maps the symmetries (\ref{39ww}) to their canonical form (\ref{cform}) is
		\begin{equation}\label{40ww}
		X= \sqrt{\frac{2 x y-1}{y^2}-1}, \quad Y = \frac{x y-1}{y}.
		\end{equation}
		This transformation also reduces Eq. (\ref{38ww}) to Eq. (\ref{e31}), with $a=1$ and $b=0$. Writing solution (\ref{28w}) in terms of $x$ and $y$ via the point transformation (\ref{40ww})  and setting  $a=1$ and $b=0$, we obtain the solution to (\ref{38ww}),
		\begin{equation}\label{}
		y= \frac{J_2}{x^2 + J_2 x + J_1 - 1}.
		\end{equation}
	}
\end{eg}

\section{Concluding remarks}\label{sec5}
There are many characterisations of linearisable second-order ODEs. One important characterisation is that every such equation is reduceable to a generic second-order ODE, Eq. (\ref{e31}), via a point transformation constructed from two suitably chosen symmetries. We have used the expanded Lie group approach to simplify Eq. (\ref{e31}) significantly via a point transformation constructed to set the parameter $b$ in Eq. (\ref{e31}) to zero. The reduced equation is subsequently mapped to the free particle equation,  $y^{\prime\prime}=0$, via another invertible point transformation  constructed by ``aligning" the respective symmetries of the two equations. The constructed point transformations are used in succession to obtain the general solution to Eq. (\ref{e31}), which is the desired proxy solution for all linearisable second-order differential equations. This allows construction of solutions of all equations in this class algorithmically using only two suitably chosen symmetries of the equation. We have illustrated the solution routine  with three examples. 

\section*{Data Availability}
No data were used to support this study.

\section*{Conflicts of Interest}

The author declares that there is no conflict of interest regarding the publication of this paper.

\section*{Acknowledgements}
The author would like to thank the Directorate of Research Development and Innovation of Walter
Sisulu University for continued financial support.

\end{document}